\documentclass{amsart}
\usepackage{amsmath,amssymb,latexsym,amsthm,amscd}
\font\goth=eusm10

\newcommand\wh{\widetilde{H}}

\newcommand\E{\mathcal E}

\newcommand\Ii{\hbox{\goth I}}

\newcommand\M{\mathcal M}

\newcommand\ZZ{\mathbf{Z}}
\newcommand\bz{\mathbf{Z}}

\newcommand\NN{\mathbf{N}}

\newcommand\Oc{\hbox{\goth O}}
\newcommand\Pj{\mathbf{P}}

\newcommand\wH{\widetilde{H}}

\newtheorem{proposition}{Proposition}
\newtheorem{theorem}{Theorem}

\newtheorem{definition}{Definition}
\newtheorem{corollary}{Corollary}
\newtheorem{remark}{Remark}

\numberwithin{proposition}{section}
\numberwithin{definition}{section}
\numberwithin{corollary}{section} \numberwithin{remark}{section}
\numberwithin{lemma}{section} \numberwithin{equation}{section}

\numberwithin{theorem}{section}

\newtheorem{question}{Question}
\numberwithin{question}{section} \numberwithin{case}{section}
\numberwithin{example}{section} \numberwithin{conjecture}{section}

\newtheorem*{dimo-p:2}{Proof of Theorem \ref{p:2}}

\begin{document}
\date{\today}
\email{carlo.madonna@uclm.es}
\address{C.G.Madonna,
Universidad de Castilla -- La Mancha,
Plaza de la Universidad 3, 02071 Albacete (Spain)}
\begin{abstract}
We give several examples of the existence of infinitely many 
divisorial conditions on the moduli space of
polarized K3 surfaces $(S,H)$ of degree $H^2=2g-2$, $g \geq 3$, 
and Picard number $\rho(S)=rk N(S)=2$, 
such that for a general K3 surface $S$ satisfying 
these conditions the moduli space of sheaves
$M_S(r,H,s)$ is birationally equivalent to the Hilbert scheme $S[g-rs]$ of zero-dimensional subschemes of $S$
of length equal to $g-rs$.
This result generalizes the main result of \cite{Nik1} when $g=rs+1$ 
and of \cite{Monat} when $r=s=2$,
$g \geq 5$.
\end{abstract}

\title[On some moduli spaces of bundles on K3 surfaces, II]{On some moduli spaces of bundles \\
on K3 surfaces, II}
\author[C.G.Madonna]{C.G.Madonna} \thanks{The author was supported by the EPSRC grant EP/D061997/1. The author is member 
of the project MTM2007-67623 founded by the Spanish MEC}
\maketitle

\section{Introduction to the main result}

Let $S \subset \Pj^g$ be a smooth projective K3 surface of genus $g \geq 3$ and $H$ be a general hyperplane section of
$S$, $H^2=2g-2$. Let us denote by $\wh(S,\ZZ)$ the Mukai lattice of $S$ and let $v=(v_0,v_1,v_2) \in \wh(S,\ZZ)$ be an algebraic primitive
Mukai vector (i.e. $\ZZ v$ is a primitive sublattice of $\wh(S,\ZZ)$ and $v_1 \in N(S)$ -- the Picard lattice of $S$).
We briefly remind that the {\it Mukai lattice} of $S$ is the full cohomology group 
\begin{equation}
\wh(S,\ZZ):=H^0(S,\ZZ) \oplus H^2(S,\ZZ) \oplus H^4(S,\ZZ) 
\end{equation}
endowed by the Mukai pairing
\begin{equation}
(v,w)=v_1 \cdot w_1-(v_0w_2+v_2w_0) \in \ZZ 
\end{equation}
where $v=(v_0,v_1,v_2), w=(w_0,w_1,w_2) \in \wh(S,\ZZ)$, 
$H^0(S,\ZZ)$ and $H^4(S,\ZZ)$ are naturally identified with $\ZZ$,
and the {\it transcendental lattice} of $S$ is defined by $T(S)=N(S)^{\perp} \subset H^2(S,\ZZ)$ i.e. 
it is the orthogonal to the
Picard lattice $N(S)$ with respect to the intersection product of $H^2(S,\ZZ)$.

We extensively studied the problem to find conditions on $S$ 
such that the moduli space $M_S(v)$ of $H$-stable sheaves over $S$ with Mukai vector $v$
is isomorphic to the K3 surface $S$ when $v$ is isotropic, i.e. $(v,v):=v^2=0$, and 
primitive.
The interested reader can refer to \cite{MN3} and the references given there for more details.
In \cite{Monat} we started to consider the case of primitive and not necessarily isotropic Mukai vectors.
Indeed, by the work of Yoshioka and Mukai (see \cite{Yos-irr}, \cite{Muk1}) the moduli space $M_S(v)$ is an irreducible symplectic
variety of dimension equal to $v^2+2$ which is equivalent under deformations to the Hilbert scheme 
$S[v^2/2+1]$ of zero-dimensional subschemes $Z$
of $S$ of length $l(Z)=v^2/2+1$. The question we want to answer here is the following:

\begin{question}
Let $v \in \wh(S,\ZZ)$ be an algebraic and primitive Mukai vector and let $M:=M_S(v)$ be the moduli space of $H$-stable sheaves
over $S$ with Mukai vector $v$. Let $M':=S[v^2/2+1]$ be the Hilbert scheme of zero-dimensional subschemes of $S$ of length $v^2/2+1$.
When is $M$ birational to $M'$?
\end{question}

In this paper we will refer to a birational map between $M$ and $M'$ as a ''birational isomorphism" or just ''isomorphism" for short.
We are mainly interested here to the case when $v^2>0$. When $v^2=0$, a complete answer was given in \cite{MN3} and
\cite{Nik1} \cite{Nik3}
for a general K3 surface $S$ of Picard number $\rho(S):=rk \ N(S) \leq 2$ (i.e. with $Aut(T(S),H^{2,0}(S))=\pm 1$).

For any $g \geq 5$ we already considered the previous question in \cite{Monat} for the case of Mukai vector
$(2,H,2)$. We found, there, infinitely many divisorial conditions on the moduli space of polarized K3 surfaces
of degree $H^2=2g-2$ for which one has a birational isomorphism $M \cong M'$.
Here we partially answer the previous question when the Mukai vector is 
\begin{equation}
v=(r,H,s)
\end{equation}
with $r,s \in \ZZ$, $r,s \geq 1$, $H^2=2g-2$,
$\rho(S)=2$, and $\gamma(H)=1$,
where $\gamma(H) \in \NN$ is defined by $H \cdot N(S) \cong \gamma(H) \ZZ$.

In this case we have that
$(v,v)=2g-2-2rs=2(g-1-rs)$
and the moduli space of $H$-stable sheaves $E$ with Mukai vector $v(E)=v=(r,H,s)$ 
\begin{equation}
M=M_S(v)
\end{equation}
has dimension equal to $2(g-rs)$. The corresponding Hilbert scheme is given by
\begin{equation}
M'=S[g-rs].
\end{equation}
We remind two beautiful examples given by Mukai in which for 
general $S$ (say $N(S) \cong \ZZ[H]$) one cannot hope to get any isomorphism
between $M$ and $M'$.
The first one corresponds to the case $g=5,r=s=2$. In this case $v$ is isotropic and the surface $S=M'$ is an octic K3 complete intersection
of three quadrics in $\Pj^5$ while $M$ is a K3 double plane, i.e. a ramified sextic double plane, 
isomorphic to the double covering of the net of quadrics whose base locus defines $S$
ramified over the locus of singular quadrics of the net (see \cite{MN1} for more details on this case).
Of course in general $M'=S \not\cong M=M_S(2,H,2)$ since $\det N(M)=2 \ne \det N(S)=8$.

The second example, also related to the geometry of quadrics, is a 4-dimensional analogue to the above one, 
corresponding to the case $g=6, r=s=2$. The interested reader can refer to \cite{Muk2} for more details.

\

In this short note we give several examples of infinitely many divisorial 
conditions on the moduli space of polarized K3 surfaces $(S,H)$ of degree
$H^2=2g-2$
for which one gets a birational isomorphism between $M$ and $M'$.
The main tool will be in considering the following isometries of the Mukai lattice of $S$: \par
-- the tensorization for a line bundle $D \in N(S)$, say, $T_D:\wh(S) \to \wh(S)$, defined by
\begin{equation}
T_D(r,c_1,s)=(r,c_1+rD,s+r(D^2/2)+D \cdot c_1)
\end{equation}
and \par
-- the reflection $\delta:\wh(S,\ZZ) \to \wh(S,\ZZ)$ defined by
\begin{equation} \label{eq:reflection}
\delta(r,c_1,s)=(s,c_1,r).
\end{equation}
Indeed when $c_1=H$ if $w:=T_D(r,H,s)$  
then one gets an isomorphism $M_S(r,H,s) \cong M_S(w)$. 
One has also an isomorphism $M_S(r,H,s) \cong M_S(s,H,r)$.

Further we will mainly need to distinguish two cases. For this we introduce the following general definition.

\begin{definition}
Let $v=(r,H,s) \in \wH(S,\ZZ)$. 
We say that $v$ is of {\it type} $+$ if there exists a divisor $D \in N(S)$ 
such that $v^+:=T_D(v)=(r,H+rD,1)$.
We say that $v$ is of {\it type} $-$ if there exists a divisor $D \in N(S)$
such that 
$v^-:=T_D(v)=(r,H+rD,-1)$.
\end{definition}

Since $M_S(v) \cong M_S(v^{\pm})$,
in order to answer our main question we will answer to the following:

\begin{question}
Under which conditions on the Picard lattice $N(S)$ of the polarized K3 surface
$S$ of degree $H^2=2g-2$ the Mukai vector $v=(r,H,s)$ is of type $\pm$ and $M_S(v^\pm) \cong S[g-rs]$?
\end{question}

If $v$ is of type $+$ then applying the reflection $\delta$,
we get that 
\begin{equation}
M_S(v) \cong M_S(r,H+rD,1)=M_S(v^+) \cong M_S(1,H+rD,r) \cong S[g-rs].
\end{equation}
Hence, in this case we are
done if we are able to find the divisor $D$ as above such that $T_D(r,H,s)=(r,H+rD,1)$. 
To find the divisor $D$ we consider a primitive element $G \in N(S)$ orthogonal
to $H$ such that $N(S)$ is generated up to its index by $H$ and $G$
(see Proposition \ref{pro:reticoli}).
Then, when $v$ is of type $-$ we 
construct explicitly the required birational isomorphism
$M_S(v^-) \cong S[g-rs]$.

\

Our first result is the following:

\begin{proposition} \label{p:2}
Let $S$ be a polarized K3 surface of degree $H^2=2g-2$, general with 
$\rho(S)=2$, $\gamma(H)=1$, $\det N(S)=-d$, $d \equiv \mu^2 \mod 4(g-1)$, $\mu \mod 2g-2 \in (\ZZ/(2g-2))^*$.
Let $v=(r,H,s)$. Assume that $d$ is not a square i.e. 
the Picard lattice $N(S)$ does not represent zero.
Then, if $d \in \mathcal D_\pm$ where
\begin{equation}
\begin{aligned}
{\mathcal D}_\pm=  \{ & d \in \NN : \ d=\frac{(rx+2(g-1))^2-4(g-1)(\pm r-rs+g-1)}{r^2 y^2}, \\
&  x,y \in \ZZ \ , \ x \equiv \mu y \mod 2g-2 \}
\end{aligned}
\end{equation}
the Mukai vector $v$ is of type $v^\pm$. Indeed 
we have that $v^\pm=(r,F,\pm 1)$ where 
\begin{equation} \label{eq:a-serie}
F:=H+rD
\end{equation}
with $F\cdot H \equiv r \mu y \mod 2g-2$ and $D=(xH+yG)/(2g-2) \in N(S)$,
where $(x,y)$ is an integral solution of the equation 
$(rx+2(g-1))^2-dr^2y^2=4(g-1)(\pm r-rs+g-1)$ with $x \equiv \mu y \mod 2g-2$.
Moreover one can choose the above divisor $D \in N(S)$ such that $x \ll 0$.
\end{proposition}

Using the above result we then derive the main result of the paper:

\begin{theorem} \label{th:main}
Let $S$ be a general K3 surface of degree $H^2=2g-2$ with Picard group $N(S):=N$ and Picard number $\rho(S)=2$,  $\gamma(H)=1$,
$\det N=-d$, $d \equiv \mu^2 \mod 4(g-1)$ and $\mu \mod 2g-2 \in (\ZZ/(2g-2))^*$.
If $d \in \mathcal D_\pm$ is not a square, 
then:
\begin{enumerate}
\item there exists a divisor $D=(xH+yG)/(2g-2) \in N$ as in previous proposition such that 
the divisor $F=H+rD \in N(S)$
satisfies
\begin{equation} \label{eq:2}
F^2=(2g-2)+r(\pm 2-2s)
\end{equation}
and $D \cdot H \ll 0$, $F \cdot H \equiv r \mu y \mod 2g-2$,
where $(x,y)$ is an integral solution of the equation
$(rx+2(g-1))^2-dr^2y^2=4(g-1)(\pm r-rs+g-1)$  with $x \equiv \mu y \mod 2g-2$;
\item there exists a birational isomorphism $M_S(r,H,s) \cong S[g-rs]$.
\end{enumerate}
\end{theorem}
 
\medskip

From the proof of the theorem it follows 
the above birational isomorphism is obtained by the composition of the isomorphism
$$
\varphi(D)_\pm:M_S(r,H,s) \to M_S(r,H+rD,\pm 1)
$$ 
induced by the tensorization 
for the line bundle $D$ as in (1) and:
\smallskip \par
\noindent
(a) when $d \in \mathcal D_+$, the
isomorphism $\delta:M_S(r,H+rD,1) \to M_S(1,H+rD,r)$ induced by the reflection $\delta$ (see (\ref{eq:reflection})), 
and the isomorphism
$\varphi_+: M_S(1,H+rD,r) \to S[g-rs]$ inverse of the isomorphism $\varphi_+^{-1}:S[g-rs] \to M_S(1,H+rD,r)$ 
defined by $\varphi_+^{-1}(Z) \cong \Ii_Z(F)$;
\smallskip\par\noindent
(b) when $d \in \mathcal D_-$, with a birational isomorphism $\varphi_-:M_S(r,H+rD,-1) \cong S[g-rs]$
inverse of a birational isomorphism $(\varphi_-)^{-1}: S[g-rs] \cong M_S(r,H+rD,-1)$ which is explicitly given in the proof of Theorem \ref{th:main}.

It is important to notice that under the conditions of the previous theorem, one can apply
the reflection $\delta$ to get
a birational isomorphism $M_S(r,H,s) \cong M_S(s,H,r)$.
Then, arguing in similar way, one can show  the following: if $d \in {\widetilde{\mathcal D}}_\pm$
is not a square, where
\begin{equation}
\begin{aligned}
{\widetilde {\mathcal D}}_\pm=\{ & d \in \NN : \ d=\frac{(sx+2(g-1))^2-4(g-1)(\pm s-rs+g-1)}{s^2 y^2}, \\
& x, y \in \ZZ \ , \ x \equiv \mu y \mod 2g-2 \}
\end{aligned}
\end{equation}
then $M_S(r,H,s) \cong S[g-rs]$. 
In particular for all such $d$ there exists a divisor $\widetilde D=(xH+yG)/(2g-2) \in N(S)$ such that 
the divisor 
\begin{equation} \label{eq:b-serie}
{\widetilde F}=H+s \widetilde D
\end{equation}
satisfies
\begin{equation} \label{eq:22}
{\widetilde F}^2=(2g-2)+s(\pm 2-2r).
\end{equation}
and ${\widetilde D} \cdot H \ll 0$, ${\widetilde F} \cdot H \equiv s \mu y \mod 2g-2$, where $(x,y)$ is an integral solution of the equation
$(sx+2(g-1))^2-ds^2y^2=4(g-1)(\pm s-rs+g-1)$  with $x \equiv \mu y \mod 2g-2$.
We leave to the reader to write a statement similar to the one of our main theorem in this case.

We remark if $r=s$ then ${\mathcal D}_{\pm}={\widetilde{\mathcal D}}_{\pm}$.

Remember also when $H^2=2g-2=2rs$, i.e. the Mukai vector $v=(r,H,s)$ is isotropic and the moduli space $M_S(v)$
is a K3 surface, we have that the corresponding divisor $F$ satisfies the conditions
\begin{equation} \label{e:isoa}
F^2=\pm 2r \ , \ F \cdot H \equiv 0 \mod r
\end{equation}
and the divisor $\widetilde F$ satisfies the conditions
\begin{equation} \label{e:isob}
{\widetilde F}^2=\pm 2s \ , \ {\widetilde F} \cdot H \equiv 0 \mod s.
\end{equation}
In this case in \cite{Nik1} 
it was proved that these two conditions (\ref{e:isoa}) and (\ref{e:isob}) 
are not only sufficient but also necessary to 
have the isomorphism $S \cong M_S(v)$ when $S$ is general with $\rho(S) \leq 2$.

Our main result also gives the following interesting:
\begin{corollary} \label{cor:infinite}
Under the conditions and notations of Theorem \ref{th:main}, let $\mathcal D=\mathcal D_+ \cup \mathcal  D_-$ and
$\widetilde{\mathcal D}=\widetilde{\mathcal D}_+ \cup \widetilde{\mathcal D}_-$. If
either $r|g-1$ or $s|g-1$ or
$r|2$ or $s|2$
then ${\mathcal D} \cup {\widetilde{\mathcal D}}$ is infinite.
\end{corollary}

Note that if $g-1=rs$, i.e. $v=(r,H,s)$ is isotropic, it then
follows that ${\mathcal D}$ and ${\widetilde{\mathcal D}}$ are infinite.
This is also the case when $g \geq 5$ and $v=(2,H,2)$ (see \cite{Monat}).

\

It is interesting to rewrite our previous result in term of Beauville-Bogomolov quadratic
form $q$ (and bilinear form $b$) on $H^2(S[g-rs],\ZZ)$
(see e.g. \cite{B} for details). 
Remember that $H^2(S[g-rs],\ZZ) \cong H^2(S,\ZZ) \oplus_\perp \ZZ f$, with $f^2=-2(g-rs-1)$.
The corollary below shows how our main result here can be considered as a natural
generalization to non isotropic Mukai vectors of the result given in \cite{Nik1} for the case of isotropic Mukai vectors. 

\begin{corollary}
Under the conditions and notations of Theorem \ref{th:main} if $d \in {\widetilde{\mathcal D}} 
\cup {\mathcal D}$
is not a square
then there exists either a divisor class $F \in N(S)$ as in (\ref{eq:a-serie}) such 
that
\begin{equation}
h_1=F+\epsilon f \in H^2(S[g-rs],\ZZ)
\end{equation}
satisfies
\begin{equation} \label{q1}
q(h_1)=\pm 2r \ , \ b(h_1,H)=F \cdot H \equiv r \mu y \mod 2g-2
\end{equation}
with $(rx+2(g-1))^2-dr^2y^2=4(g-1)(\pm r-rs+g-1)$ and $x \equiv \mu y \mod 2g-2$,
or a divisor class ${\widetilde F}$ as in (\ref{eq:b-serie}) such that the element 
\begin{equation}
{\widetilde h}_1={\widetilde F}+\epsilon f \in H^2(S[g-rs],\ZZ)
\end{equation}
satisfies
\begin{equation} \label{q2}
q({\widetilde h}_1)=\pm 2s \ , \ b({\widetilde h}_1,H)={\widetilde F} \cdot H \equiv s \mu y \mod 2g-2
\end{equation} 
with $(sx+2(g-1))^2-ds^2y^2=4(g-1)(\pm s-rs+g-1)$ and 
$x \equiv \mu y \mod 2g-2$.. Here
$\epsilon=0$ when $g=rs+1$ and $\epsilon=1$ when $g>rs+1$.
\end{corollary}

{\it Acknowledgement.} This paper was written during the author's visit to the Mathematics
Department of The University of Liverpool. The author thanks Professor Viacheslav Nikulin for 
useful discussions.

\section{A proof of the main result}

In previous section we stated our main result. 
We also made several remarks, derived some corollaries, and made some
questions we hope to answer in further works. It remains now to give a proof of Theorem \ref{th:main}.
In this section we will do this.
We will now give the proof of Proposition \ref{p:2} which will be needed
for the proof of Theorem \ref{th:main}. 

At first, to fix notations and terminology, we will recall some
basic facts on lattices (see \cite{Nik2} for more details) we will need in the sequel.
A {\it lattice} $L$ is a
non-degenerate integral symmetric bilinear form. I. e. $L$ is a
free $\ZZ$-module equipped with a symmetric pairing $x\cdot y\in
\bz$ for $x,\,y\in L$, and this pairing should be non-degenerate.
We denote $x^2=x\cdot x$. 
The {\it determinant} of
$L$ is defined to be $\det L=\det(e_i\cdot e_j)$ where $\{e_i\}$
is some basis of $L$. We will use notation $[x_1,...,x_n]$ to mean the lattice generated by
the elements $x_1,...,x_n$. 

An {\it isometry} of lattices is an isomorphism of modules which preserves the
product.

Let $S$ be a polarized K3 surface of degree $H^2=2g-2$.
For the rest of the paper we will assume that the following conditions on the Picard lattice
$N(S)$ of $S$ are satisfied:
\begin{equation} \label{e:*}
\begin{cases}
\gamma(H)=1 \\
\rho(S) = 2 \\
\mu \mod 2g-2 \in (\ZZ/(2g-2))^* \\
\ZZ v \subset \wh(S,\ZZ) \ \text{is primitive} 
\end{cases}
\end{equation}

We also need the following (see \cite{Nik1} for a proof):

\begin{proposition} \label{pro:reticoli}
Let $S$ be a K3 surface with a polarization $H$ of degree $H^2=2g-2$. Assume that conditions 
(\ref{e:*}) hold. Then the lattice $N(S)$ is defined by its determinant $\det N(S)=-d$, where $d \in \NN$,
and $\mu^2 \equiv d \mod 4(g-1)$. There exists a unique choice of a primitive vector $G$ orthogonal to $H$ such that
$$
N(S)=[H,G, (\mu H+G)/(2g-2)]
$$
where $G^2=-(2g-2)d$, and $\mu \mod 2g-2 \in (\ZZ / (2g-2))^*$.
It follows that 
$$
N(S)=\{ (xH+yG)/(2g-2) : x,y \in \ZZ, x \equiv \mu y \mod 2g-2 \}.
$$
\end{proposition}

We can now derive the following:

\begin{proposition} 
Let $S$ be a K3 surface with a polarization $H$ of degree $H^2=2g-2$. Suppose that conditions 
(\ref{e:*}) hold. Let $d \in \mathcal D_\pm$, where 
$$
\begin{aligned}
\mathcal D_-=\{ & d \in \NN : \ d=\frac{(rx+2(g-1))^2-4(g-1)(-r-rs+g-1)}{r^2y^2}, \\
& x,y \in \ZZ, x \equiv \mu y \mod 2g-2 \}
\end{aligned}
$$
$$
\begin{aligned}
\mathcal D_+=\{ & d \in \NN : \ d=\frac{(rx+2(g-1))^2-4(g-1)(r-rs+g-1)}{r^2y^2}, \\
& x,y \in \ZZ, x \equiv \mu y \mod 2g-2 \}.
\end{aligned}
$$
Moreover suppose that $d$ is not a square. Then there exists a 
divisor $D=(xH+yG)/(2g-2) \in N(S)$ such that 
$T_D(r,H,s)=(r,H+rD,\pm 1)$ and $D \cdot H \ll 0$,
where $(x,y)$ is an integral solution of the equation 
$(rx+2(g-1))^2-dr^2y^2=4(g-1)(\pm r-rs+g-1)$ with $x \equiv \mu y \mod 2g-2$.
For all such $d$ the divisor $D$ as above is then given by $D=(xH+yG)/(2g-2)$ 
and  we have that $F:=H+rD$ satisfies
$$
F^2=(2g-2)+r(\pm 2-2s)
$$
with $F\cdot H \equiv r \mu y \mod 2g-2$.
\end{proposition} 

\begin{proof}
By the previous proposition we can write a divisor $D \in N(S)$ in the following form
$$
D=\frac{xH+yG}{2g-2}
$$
with $x,y \in \ZZ$ and $x \equiv \mu y \mod 2g-2$.
Then $T_D(r,H,s)=(r,H+rD,\pm 1)$ if and only if 
$(x,y)$ is an integral solution of the Pell--type equation
\begin{equation} \label{eq:pell}
(rx+2(g-1))^2-d(ry)^2=4(g-1)(\pm r-rs+g-1)
\end{equation}
In particular if we denote $H+rD$ by $F$ we have that
$F^2=(2g-2)+r(\pm 2-2s)$ and $F \cdot H \equiv r \mu y \mod 2g-2$.
Finally, since the solutions of the above equation are not bounded when $d$ is not a square, then
we may assume that $D \cdot H \ll 0$.
\end{proof}

We can now give a proof of our main result

\begin{theorem} \label{th:main-i}
Let $S$ be a general K3 surface of degree $H^2=2g-2$ with Picard group $N(S):=N$ and Picard number $\rho(S)=2$,  $\gamma(H)=1$,
$\det N=-d$, $d \equiv \mu^2 \mod 4(g-1)$ and $\mu \mod 2g-2 \in (\ZZ/(2g-2))^*$.
If $d \in \mathcal D_\pm$ is not a square, 
then:
\begin{enumerate}
\item there exists a divisor $D=(xH+yG)/(2g-2) \in N$ as in previous proposition such that 
the divisor $F=H+rD \in N(S)$
satisfies
\begin{equation} \label{eq:2i}
F^2=(2g-2)+r(\pm 2-2s)
\end{equation}
and $D \cdot H \ll 0$, $F \cdot H \equiv r \mu y \mod 2g-2$,
where $(x,y)$ is an integral solution of the equation
$(rx+2(g-1))^2-dr^2y^2=4(g-1)(\pm r-rs+g-1)$  with $x \equiv \mu y \mod 2g-2$.;
\item there exists a birational isomorphism $M_S(r,H,s) \cong S[g-rs]$.
\end{enumerate}
\end{theorem}

\begin{proof}
\normalfont{
We already noticed that if $d \in \mathcal D_+$ then the result follows soon.
Indeed in this case the birational isomorphism $M' \cong M$ is obtained as the composition 
$\varphi_+ \circ \delta \circ \varphi(D)_+$
of
the isomorphism
\begin{equation}
\varphi(D)_+: M_S(r,H,s) \cong M_S(r,H+rD,1)
\end{equation}
given by the tensorization for the line bundles $D \in N=N(S)$ of previous proposition, the isomorphism
\begin{equation}
\delta: M_S(r,H+rD,1) \cong
M_S(1,H+rD,r) 
\end{equation}
induced by the reflection (\ref{eq:reflection}),
and the isomorphism 
\begin{equation}
\varphi_+:M_S(1,H+rD,r) \cong S[g-rs].
\end{equation}
which is the inverse of the isomorphism $\varphi_+^{-1}:S[g-rs] \to M_S(1,H+rD,r)$ defined by
$\varphi_+^{-1}(Z)=\Ii_Z(H+rD)$.
Therefore we assume that $d \in \mathcal D_-$.
Let $v=(r,H,s)$ be of type $-$ and hence let $v^-=(r,F,-1)$ where 
$F=H+rD$ with $D \in N(S)$ as in the previous proposition, in particular such that $D \cdot H \ll 0$.
We will give here explicitly the map $(\varphi_-)^{-1}: S[g-rs] \to M_S(v^-)$.
Then we will get the required birational isomorphism as composition of 
the isomorphism
\begin{equation}
\varphi(D)_-: M_S(r,H,s) \cong M_S(r,H+rD,-1)
\end{equation}
given by the tensorization for the line bundles $D \in N(S)$ as in previous proposition and the birational isomorphism
\begin{equation}
\varphi_-:M_S(r,H+rD,-1) \cong S[g-rs].
\end{equation}
Let $Z \in S[g-rs]$ be a general subscheme of $S$ of length $l(Z)=g-rs$.
Then $h^1 \Ii_Z(F)=r-1$.
Indeed we have the short exact sequence defining $Z \subset S$
$$
0 \to \Ii_Z(F) \to \Oc_S(F) \to \Oc_Z(F) \to 0
$$
and the corresponding exact sequence in cohomology.
It then follows by standard calculations that $h^1 \Ii_Z(F) \geq r-1$.
Following the construction of \cite{Tyurin1}
we then have the exact sequence 
$$
0 \to H^1 \Ii_Z(F) \otimes \Oc_S \to \E=\E_{Z} \to \Ii_Z(H+rD) \to 0
$$
where $\E$ has $rk(\E)=h^1 \Ii_Z(F)+1$, $c_1(\E)=F:=H+rD$ and
$\chi(\E)=r-1$, $h^1 \E=0$. Hence the Mukai vector of $\E$ is $v(\E)=(r,F,-1)$.
We define here $(\varphi_-)^{-1}(Z)=\E_Z=\E$.
It remains to show that $\E$ is $H$-stable.
We will show that $\E(-D)$ is $H$-stable.
Remember  that $F=H+rD$ for a suitable $D \in N(S)$ such
that $D \cdot H \ll 0$. 
We then tensorize the above exact sequence by $\Oc_S(-D)$ to get
$$
0 \to \Oc_S^{r-1}(-D) \to \E(-D) \to \Ii_Z(H+(r-1)D) \to 0
$$
Now let us assume that $\E(-D)$ is unstable.  Then consider the maximal
destabilizing subsheaf of $\E(-D)$, say $\M$.
We have:
$$
\mu(\M)=\frac{(c_1(\M) \cdot H)}{rk(\M)} > \mu(\E(-D))=\frac{2g-2}{r}>0.
$$
Let $i$ be the inclusion map of $\M$ in $\E(-D)$ and let $\beta:\E(-D) \to
\Ii_Z(H+(r-1)D)$ be the map given in
the above exact sequence.
Then the image of $\M$ in $\Ii_Z(H+(r-1)D)$ by the composite map $\beta
\circ i$ is a subsheaf of rank 1 which is
of type $\Ii_W(H+(r-1)D)$ for some zero-dimensional
subscheme $W$ containing $Z$.
Then we have an exact sequence
$$
0 \to \M' \to \M \to \Ii_W(H+(r-1)D) \to 0
$$
and $c_1(\M')=c_1(\M)-(H+(r-1)D)$.
We then have
\begin{equation}
\begin{aligned}
\mu(\M') & =\frac{c_1(\M') \cdot H}{rk(\M')}
=\frac{c_1(\M) \cdot H-(H+(r-1)D) \cdot H}{rk(\M')}> \\
& > \frac{c_1(\M) \cdot H}{rk(\M')}
\geq \frac{c_1(\M)\cdot H}{rk(\M)}=\mu(\M)
\end{aligned}
\end{equation}
since $c_1(\M) \cdot H>0$, $D \cdot H \ll 0$, $rk(\M')=rk(\M)-1$. Hence $\mu(\M')>\mu(\M)$.
This is impossible since we assumed $\M$ to be maximal destabilizing.}
\end{proof}

From the proof of previous theorem 
we have the following more general statement:

\begin{theorem} \label{th:main2}
Let $S$ be a K3 surface of degree $H^2=2g-2$ and suppose that $N(S)$ contains a
rank 2 primitive sublattice $N \subset N(S)$ with $H \in N$. Assume that 
$(S,N)$ satisfies the conditions of Theorem \ref{th:main}. With same notations of Theorem \ref{th:main} for $\mu, d, \gamma, \mathcal D_{\pm}$,
if $d \in \mathcal D_\pm$ is not a square, 
then:
\begin{enumerate}
\item there exists a divisor $D=(xH+yG)/(2g-2) \in N$ such that 
the divisor $F=H+rD$
satisfies
\begin{equation} \label{eq:3}
F^2=(2g-2)+r(\pm 2-2s)
\end{equation}
and $D \cdot H \ll 0$, $F \cdot H \equiv r \mu y \mod 2g-2$, where 
$(x,y)$ is an integral solution of the equation
$(rx+2(g-1))^2-dr^2y^2=4(g-1)(\pm r-rs+g-1)$  with $x \equiv \mu y \mod 2g-2$;
\item there exists a birational isomorphism $M_S(r,H,s) \cong S[g-rs]$.
\end{enumerate}
\end{theorem}

\begin{remark}
\normalfont{
When $v^2=0$ Nikulin gave in \cite{Nik3} examples of K3 surfaces $S$ with Picard group $N(S)$ of rank $\rho(S)=3$ such
that $M_S(v) \cong S$ while $N(S)$ does not contain any rank 2 primitive sublattice $N$ as in previous theorem.
In other words, examples when $\rho=3$ in which one has an isomorphism which 
is not a specialization of an isomorphism between K3 surfaces of Picard number $2$.
These examples have all $\gamma(H) \neq 1$.
We do not know if one can find similar examples when $\gamma(H)=1$ (with isotropic or non isotropic Mukai vectors as well).
}
\end{remark}

\section{Examples}

Let us consider the cases in which $v^2=0$ i.e. $g=5$ and $v=(2,H,2)$.
Our main result in this case reads as follows:

\begin{theorem} 
Let $S$ be a general K3 surface of degree $H^2=8$ with Picard group $N(S):=N$ and Picard number $\rho(S)=2$,  $\gamma(H)=1$,
$\det N=-d$, $d \equiv \mu^2 \mod 16$ and $\mu \mod 8 \in (\ZZ/8)^*$, i.e. $\mu=1$ if $d \equiv 1 \mod 16$ and $\mu=3$ if $d \equiv 9 \mod 16$
(see Section 3.1 of \cite{MN1}).
If $d \in \mathcal D_\pm$ is not a square, 
then:
\begin{enumerate}
\item there exists a divisor $D=(xH+yG)/8 \in N$ such that 
the divisor $F=H+2D$
satisfies
\begin{equation} \label{eq:2ii}
F^2=\pm 4
\end{equation}
and $D \cdot H \ll 0$, $F \cdot H \equiv 2 \mu y \mod 8$,
where $(x,y)$ is an integral solution of the equation
$(x+4)^2-dy^2=\pm 8$ with $x \equiv \mu y \mod 8$;
\item there exists a birational isomorphism $M_S(2,H,2) \cong S$.
\end{enumerate}
\end{theorem}

In this cases the the families $\mathcal D$ and $\mathcal{\tilde{D}}$ coincide, moreover by Corollary \ref{cor:infinite}
$\mathcal D$ is infinite.
The two families $\mathcal D_+$ and $\mathcal D_-$ are given by
$$
\mathcal D_\pm=\{ d \in \NN : \ d=\frac{(x+4)^2 \mp 8}{y^2}, 
x,y \in \ZZ, x \equiv \mu y \mod 8 \}
$$
After changing $x+4=a$ and $y=b$ then we get the following equations
$$
a^2-db^2=\pm 8
$$
with $a-4 \equiv \pm b \mod 8$ if 
$d \equiv 1 \mod 16$ and $a-4 \equiv \pm 3 b \mod 8$ if 
$d \equiv 9 \mod 16$. When $b=1$ we then get the conditions 
$a \equiv \pm 3 \mod 8$ if $d \equiv 1 \mod 16$ and $a \equiv \pm \mod 8$ if
$d \equiv 9 \mod 16$.

We derive the following possible values for $d$
$$
17,33,41,57,89,73,113,129,161,177, ...
$$
which is an infinite subfamily we already found in \cite{MN1}.

\smallskip

When $v=(2,H,2)$ and $H^2 \geq 10$ then $v^2>0$ is not anymore 
isotropic.
In this case from our main result we derive the following theorem
already showed in \cite{Monat}

\begin{theorem} 
Let $S$ be a general K3 surface of degree $H^2=2g-2 \geq 10$ 
with Picard group $N(S):=N$ and Picard number $\rho(S)=2$,  $\gamma(H)=1$,
$\det N=-d$, $d \equiv \mu^2 \mod 4(g-1)$ and 
$\mu \mod 2g-2 \in (\ZZ/(2g-2))^*$.
If $d \in \mathcal D_\pm$ is not a square, 
then:
\begin{enumerate}
\item there exists a divisor $D=(xH+yG)/(2g-2) \in N$ such that 
the divisor $F=H+2D$
satisfies
\begin{equation} \label{eq:2iii}
F^2=2(g-7)
\end{equation}
or
\begin{equation}
F^2=2(g-3)
\end{equation}
and $D \cdot H \ll 0$, $F \cdot H \equiv 2 \mu y \mod 2g-2$, where
$(x,y)$ is an integral solution of the equation
$(x+g-1)^2-dy^2=(g-1)(g-5 \pm 2)$;
\item there exists a birational isomorphism $M_S(2,H,2) \cong S[g-4]$.
\end{enumerate}
\end{theorem}


\begin{thebibliography}{ADSE}

\bibitem{B} A.Beauville, \emph{Vari\'et\'es K\"ahleriennes dont la premi\`ere classe de Chern est nulle}, J. Differential Geom. \textbf{18}(1983), no. 4, 755--782.

\bibitem{Monat} C.G. Madonna,
\emph{On some moduli spaces of bundles on K3 surfaces}, Monatsh.
Math. \textbf{146}(2005), 333--339.

\bibitem{MN1} C. Madonna and V.V. Nikulin, \emph{On a classical
correspondence between K3 surfaces}, Proc. Steklov Inst. of Math.
\textbf{241}(2003), 120--153.

\bibitem{MN3} C.G. Madonna and V.V. Nikulin,
\emph{Explicit correspondences of a K3 surface with itself}, 
Izv. Math. \textbf{72}(2008), no.3,  497--508.

\bibitem{Muk1} Sh. Mukai, \emph{On the moduli space of bundles on K3 surfaces I},
in: Vector bundles on algebraic varieties,
Tata Inst. Fund. Res. Studies in Math. no. \textbf{11} (1987),
341--413.

\bibitem{Muk2} Sh. Mukai
\emph{Moduli of vector bundles on $K3$ surfaces and symplectic manifolds}, Sugaku Expositions \textbf{1}(1988), no. 2, 139--174.

\bibitem{Nik2} V.V. Nikulin, 
\emph{Integral symmetric bilinear forms and some of
their geometric applications}, 
Math USSR-Izv. \textbf{14}(1980), no. 1, 103--167.

\bibitem{Nik1} V.V. Nikulin, \emph{On correspondences of a K3 surface
with itself. I}, Proc. Steklov Inst. of Math. \textbf{246} (2004),
204--226.

\bibitem{Nik3} V.V. Nikulin, \emph{On Correspondences of a K3 Surface with itself, II}, 
Contemporary Mathematics
\textbf{422}(2007), 121-172.

\bibitem{Nik4} V.V.Nikulin, \emph{Self-correspondences of a K3 surface via moduli of sheaves}, 
In: Y. Tschinkel and Y. Zarhin (eds.)  Algebra, Arithmetic, and Geometry. Volume II: In Honor of Yu. 
I. Manin, 2009, XII, 704 pp.,  Birkhauser Boston book, 439--464.

\bibitem{Tyurin1} A.N. Tyurin,
\emph{Cycles, curves and vector bundles on algebraic surfaces},
Duke Math. J. \textbf{54} (1987), no. 1,
1--26.

\bibitem{Yos-irr} K.Yoshioka, \emph{Some examples of Mukai's reflections on $K3$ surfaces}, J. Rein Angew. Math.  \textbf{515}(1999), 97--123.


\


\end{thebibliography}
\end{document}